\documentclass[11pt]{article}

\def\BC#1#2{\left({#1\atop#2}\right)}  
\def\abs#1{\ensuremath\left|#1\right|} 

\def\eqref#1{\mbox{(\ref{eq:#1})}}    
\def\Prob{{\ensuremath{\sf Prob}\thinspace}}
\def\Exp{{\ensuremath{\sf E}\thinspace}}

\date{November 23, 1999}

\begin{document}

\title{Walking into an Absolute Sum  \vspace{-0.25cm}}
\author{Hans J. H. Tuenter \\   
       {\footnotesize  Schulich School of Business, York University,}   \\
       {\footnotesize  Toronto, Ontario, Canada M3J 1P3}    \\
       {\footnotesize  email: htuenter@schulich.yorku.ca}}
\maketitle

\section{Introduction}
Recently, it was asked by Paul Bruckman~\cite{Bruckman:1999} to show that the sum
\begin{equation}
  S_r(n)=\sum_{k=0}^{2n} \BC{2n}{k}\abs{n-k}^r
  \label{eq:1}
\end{equation}
evaluates to $n^2\BC{2n}{n}$ for $r=3$. 
In the published solution~\cite{Strazdins:2000}, it was also noted that $S_1(n)=n\BC{2n}{n}$,
and, as a consequence, it was conjectured that $S_{2r+1}(n)$ equals the product of $\BC{2n}{n}$
and a monic polynomial of degree~$r+1$.

We show this conjecture to be true, albeit with the modification of discarding the adjectival modifier 
``monic.'' 
In fact, we show that 
\[ S_{2r+1}(n)=P_r(n)\,n\BC{2n}{n} \ \ {\rm and} \ \ S_{2r}(n)=Q_{r}(n)\,2^{2n-r},  \]
where $P_{r}(n)$ and $Q_r(n)$ are both polynomials of degree~$r$ with integer coefficients.
We then investigate the relationship of these polynomials to the 
Dumont-Foata polynomials~\cite{DumontFoata:1976}. These are generalizations of the Gandhi polynomials, which find their origin in a representation of the 
Genocchi numbers, first conjectured by Gandhi~\cite{Gandhi:1970}.
Finally, we show that the sums $S_r(n)$ are essentially the moments of a random variate,
measuring the absolute distance to the origin in a symmetric Bernoulli random walk,
after $2n$ time steps.

\section{Derivation}
We note that the sum can be rewritten as
\[ S_r(n)=2\sum_{k=0}^{n}\BC{2n}{n-k}k^r - \BC{2n}{n}\delta_{r0},
\]
with $\delta_{r0}$ the Kronecker delta.
Now consider, for $r\ge1$,
\[ n^2S_r(n)-S_{r+2}(n)
   = 2\sum_{k=0}^{n-1}\BC{2n}{n-k}k^r\left(n^2-k^2\right) 
   =  4n(2n-1)\sum_{k=0}^{n-1}\BC{2n-2}{n-1-k}k^r,
\]
leading directly to the recursion 
\begin{equation}
 S_{r+2}(n)= n^2S_r(n) - 2n(2n-1)S_r(n-1). 
 \label{eq:2}
\end{equation}
For $r=0$ the derivation is slightly more elaborate, because we need to keep track of the additional term, but leads to the same recursion,
so that~\eqref{2} is valid for all nonnegative integers~$r$.
To start the recursion, we find the value $S_0(n)=2^{2n}$ by an application of the
binomial theorem to~\eqref{1}. 
The value of $S_1(n)$ is easily obtained by breaking up the summand~$k$ to create two sums:
\[ S_1(n) =   \sum_{k=0}^n\BC{2n}{n-k}\left[\vphantom{(^0}(n+k)-(n-k)\right]  
          = 2n\sum_{k=0}^n\BC{2n-1}{n-k}-2n\sum_{k=0}^{n-1}\BC{2n-1}{n-k-1}, 
\]
and one sees that, after changing the range of summation of the second sum to start at $k=1$, all terms cancel out, 
with the exception of the summand $2n\BC{2n-1}{n}$.  
Rearranging terms gives the desired $S_1(n)=n\BC{2n}{n}$.

It is now clear that the structure of the sum depends upon the parity of $r$.
Starting with the odd values, we simplify the recursion~\eqref{2} by the substitution
$S_{2r+1}(n)=P_{r}(n)\,n\BC{2n}{n}$ to give
\begin{eqnarray}
  P_{r+1}(n) &=& n^2\left[\vphantom{(^0} P_{r}(n)-P_{r}(n-1)\right]
   +nP_{r}(n-1),
 \label{eq:3a}
\end{eqnarray}
with initial condition $P_0(n)=1$. An inductive argument now shows that $P_{r}(n)$ is a polynomial of degree~$r$ with integer coefficients,
and proves the modified conjecture.
It is not difficult to show that $r!$ is the leading coefficient of $P_{r}(n)$,  
and, hence, that these polynomials are not monic.
In fact, the only cases for which the leading coefficient is~$1$ are $r=0$ and $r=1$.
The first few polynomials are now easily determined as:
\begin{eqnarray*}
  P_0(n)&=&1, \\
  P_1(n)&=&n,  \\
  P_2(n)&=&(2n-1)n,   \\
  P_3(n)&=&(6n^2-8n+3)n,   \\
  P_4(n)&=&(24n^3-60n^2+54n-17)n,   \\
  P_5(n)&=&(120n^4-480n^3+762n^2-556n+155)n.   
\end{eqnarray*}
For the even sums we substitute $S_{2r}(n)=Q_{r}(n)\,2^{2n-r}$
to give the recursion
\begin{eqnarray}
  Q_{r+1}(n)&=&2n^2\left[\vphantom{(^0}Q_{r}(n)-Q_{r}(n-1)\right]+nQ_{r}(n-1),
  \label{eq:4a}
\end{eqnarray}
with initial condition~$Q_0(n)=1$.
This shows that $Q_{r}(n)$ is a polynomial of degree~$r$ with integer coefficients. 
It is not difficult to establish that the leading coefficient is given 
by $(2r-1)\cdot(2r-3)\cdots3\cdot1=(2r)!/(2^rr!)$,
and, hence, that these polynomials are also not monic.
Applying the recursion gives the first few polynomials as:
\begin{eqnarray*}
  Q_0(n)&=&1, \\
  Q_1(n)&=&n,\\
  Q_2(n)&=&(3n-1)n, \\
  Q_3(n)&=&(15n^2-15n+4)n, \\
  Q_4(n)&=&(105n^3-210n^2+147n-34)n, \\
  Q_5(n)&=&(945n^4-3150n^3+4095n^2-2370n+496)n. 
\end{eqnarray*}
It is worth noting that, by evaluating $S_r(n)$ for particular values of $n$,
one can derive various properties of [the coefficients of] the polynomials $P_r(n)$ and $Q_r(n)$.
For instance, it is not difficult to show that the coefficients of $P_{r}(n)$ sum to unity, 
and those of $Q_r(n)$ to $2^{r-1}$ (for $r\ge1$),
by evaluating the sums for $n=1$.
Indeed, one can derive the closed-form solutions for $S_{2r}(n)$ and $S_{2r+1}(n)$,
by solving a system of linear equations in $r$ unknowns,
representing the coefficients of the corresponding polynomial.

In the constant of the polynomials $P_r(n)/n$ 
one recognizes 
the Genocchi numbers~\cite{Comtet:1974,GrahamKP:1989}, 
named after the Italian mathematician Angelo Genocchi (1817--1889):
\[ G_2=-1,\ \ G_4=1,\ \ G_6=-3,\ \ G_8=17,\ \ G_{10}=-155,\ \ G_{12}=2073,\ \ \ldots. \]
These integers are defined through the exponential generating function
\[ \frac{2t}{e^t+1}=t+\sum_{r\ge1} G_{2r}\frac{t^{2r}}{(2r)!},
\]
and are related to the Bernoulli numbers by $G_{2r}=2(1-2^{2r})B_{2r}$.  
The Genocchi numbers are listed as sequence A001469 in the on-line version of 
the encyclopedia of integer sequences~\cite{Sloane:2000},
where additional references may be found.
The constant of the polynomials $Q_r(n)/n$ matches the first terms of the 
sequence A002105 in~\cite{Sloane:2000}, which is generated by $2^{r-1}G_{2r}/r$,
and related to the tangent numbers.
The connection to the Genocchi numbers 
will be further explored in the next section, 
where the polynomials $P_r(n)$ and $Q_r(n)$ are found to be related to special cases of the Dumont-Foata polynomials. 

Another matter of interest is the leading coefficient of the polynomials,
characterizing the behavior of the sums $S_r(n)$ for large values of $n$.
For the even indexed sums, this is easily established as 
\begin{equation}
  S_{2r}(n)
   \sim  \frac{(2r)!}{2^{2r}r!}\,2^{2n}\,n^r,
  \label{eq:3}
\end{equation}
and for the odd indexed sums we can use Stirling's formula to give $\BC{2n}{n}\sim2^{2n}/\sqrt{\pi n}$,
so that
\begin{equation} 
  S_{2r+1}(n)  
   \sim  \frac{r!}{\sqrt{\pi}}\,2^{2n}\,n^{r+\frac12}. 
  \label{eq:4}
\end{equation}
In these expressions, one can recognize the moments of a central chi-distribution,
see for instance~\cite[pp.~420--421]{JohnsonKB:1994}.
That this is no coincidence will be shown in section~\ref{Sect:RandomWalk} where we establish
the connection between the sums $S_r(n)$ and the distance to the origin 
in a symmetric Bernoulli random walk.

\section{Dumont-Foata polynomials}
In this section we show that the polynomials $P_r(n)$ and $Q_r(n)$ are related
to special cases of the Dumont-Foata polynomials~\cite{DumontFoata:1976}.
These are defined recursively by means of 
\[  F_{r+1}(x,y,z)=(x+z)(y+z)F_r(x,y,z+1)-z^2F_r(x,y,z),
\]
with initial condition $F_1(x,y,z)=1$.
Explicit expressions for these polynomials and their generating functions have been derived by Carlitz~\cite{Carlitz:1980}, but are too lengthy to display here.

The Dumont-Foata polynomials can be regarded as generalizations of the Gandhi polynomials, 
see for instance~\cite{Dumont:1974,Strehl:1979}, which are defined by
the recursion
\[ \tilde P_{r+1}(z)=(z+1)^2\tilde P_r(z+1)-z^2\tilde P_r(z),
\]
with initial condition $\tilde P_1(z)=1$.
\begin{table}[t]
\begin{center}\footnotesize
\begin{tabular}{ccccccccccc}
&&&&&    1 \\
&&&&   2 && 1 \\
&&& 6 && 8 && 3 \\
&& 24 &&  60 && 54 && 17 \\
& 120&&480&&762&&556&&155  \\
  720&&4200&&10248&&12840&&8146&&2073                                            
\end{tabular}
\end{center}
\caption{\small
  Coefficients of the Gandhi polynomials, arranged in triangular form.}
\label{Tab:Gandhi}
\end{table}
The coefficients of the first few of these polynomials are displayed in Table~\ref{Tab:Gandhi},
and can also be found in~\cite[Seq.~A036970]{Sloane:2000}.
The Gandhi polynomials arose from a conjecture made by Gandhi~\cite{Gandhi:1970}, concerning a representation of the Genocchi numbers. 
Gandhi's conjecture that $\tilde P_{r}(0)=(-1)^rG_{2r}$ was proved by 
Carlitz~\cite{Carlitz:1972}, and also by Riordan and Stein~\cite{RiordanStein:1973}.
Another polynomial that can be derived as a special case of the Dumont-Foata polynomials is given by 
$\tilde Q_r(z)=2^{r-1}F_r(\frac12,1,z)$, and is generated by the recursion
\[ \tilde Q_{r+1}(z)=(2z+1)(z+1)\tilde Q_r(z+1)-2z^2\tilde Q_r(z),
\]
with initial condition $\tilde Q_1(z)=1$.
The coefficients of the first few of these polynomials are displayed in Table~\ref{Tab:HTpol}.
For $r\ge1$, one can easily verify by substitution in~\eqref{3a} and~\eqref{4a} 
that $P_r(n)=(-1)^{r-1}n\tilde P_r(-n)$
and $Q_r(n)=(-1)^{r-1}n\tilde Q_r(-n)$. 
This gives the connection to the Dumont-Foata polynomials
(for positive $r$) as:
\begin{eqnarray*}
  P_r(n)&=&(-1)^{r-1}n\,F_r(1,1,-n)  \\ 
\noalign{\noindent and}
  Q_r(n)&=&(-2)^{r-1}n\,F_r(\textstyle\frac12,1,-n).
\end{eqnarray*}
The occurrence of the Genocchi numbers in the constant of the polynomials $P_r(n)/n$ 
is now seen to be a direct consequence of Gandhi's conjecture that $F_r(1,1,0)=(-1)^rG_{2r}$.
The occurrence of the Genocchi numbers in the constant of the polynomials $Q_r(n)/n$ 
is conjectured by the present author in the form $F_r(\frac12,1,0)=(-1)^r G_{2r}/r$.
\begin{table}[t]
\begin{center}\footnotesize
\begin{tabular}{ccccccccccc}
&&&&&    1 \\
&&&&   3 && 1 \\
&&& 15 && 15 && 4 \\
&& 105 && 210 && 147 && 34\\
&945 && 3150 && 4095 && 2370 && 496 \\
10395 && 51975 && 107415 && 111705 && 56958 && 11056 
\end{tabular}
\end{center}
\caption{\small Coefficients of the polynomials $\tilde Q_r(z)$, arranged in triangular form.}
\label{Tab:HTpol}
\end{table}

\section{Symmetric Bernoulli random walks}  \label{Sect:RandomWalk}
In a symmetric Bernoulli random walk, one considers the movements of a particle starting at time $t=0$ at the origin.
Its movements are determined by a chance mechanism, where a fair coin is flipped and 
the particle is moved one unit to the right if it is heads up, 
and one unit to the left if it is tails up. 
A more exhaustive description and in-depth study of random walks can be found 
in Feller~\cite{Feller:1950} or R\'ev\'esz~\cite{Revesz:1990}.
A more playful introduction to the topic is given in the monograph by Dynkin and Uspenskii~\cite{DynkinUspenskii:1963}.
A topic of interest is the position of the particle after $2n$ coin tosses: $Y_{2n}=X_1+X_2+\cdots+X_{2n}$,
where $X_i$ is $+1$ or $-1$ depending upon whether or not the coin showed heads in the $i$th coin toss.
Note that the $X_i$ are independent and identically distributed variates with mean~$0$ and variance~$1$.
The probability distribution of the position of the particle after $2n$ moves can 
be derived from a simple combinatorial argument, 
see for instance~\cite[p.~75]{Feller:1950} or~\cite[p.~13]{Revesz:1990},
and is given by
\[ \Prob(Y_{2n}=2k) = \BC{2n}{n-k}2^{-2n},
\]
where $k=-n,-n+1,\ldots,n$, and $n$ a positive integer.
The matter of interest in the context of this note is the distance to the origin,
$\abs{Y_{2n}}$ at time $t=2n$.
Its moments are given by
\[ \Exp{\abs{Y_{2n}}^r}=\sum_{k=-n}^{n} \BC{2n}{n-k}2^{-2n}\abs{2k}^r,
\]
and one sees that $\Exp{\abs{Y_{2n}}^r}=2^{r-2n}S_r(n)$, thus establishing the connection to the absolute sums
from the introduction.
The limit behavior of these sums now becomes clear.
By the central limit theorem, see for instance~\cite[p.~18]{Gut:1988}, 
one has that $Y_{2n}$, for sufficiently large~$n$, 
follows a normal distribution with mean $0$ and
variance $2n$.
This implies that asymptotically, $\abs{Y_{2n}}$ has a half-normal or central chi-distribution,
so that 
\[ \Exp{\abs{Y_{2n}}^r}\sim \frac{\Gamma\left[(r+1)/2\right]}{\Gamma(1/2)}\, 2^r n^{r/2}, \]
see for instance~\cite[pp.~420--421]{JohnsonKB:1994}.
This gives the asymptotic behavior of the sums as 
\[  S_r(n)=2^{2n-r}\Exp{\abs{Y_{2n}}^r}
   \sim\frac{\Gamma\left[(r+1)/2\right]}{\Gamma(1/2)}\, 2^{2n} n^{r/2}, \]
and upon expanding the gamma functions one recovers the limit results~\eqref{3} and~\eqref{4}.

\section{Discussion}
One could possibly use the relation of the Gandhi polynomials to the sums $S_{2r+1}(n)$ to gain new insight into the former.
In particular, one now has an expression to derive the function values 
of the Gandhi polynomials for negative, integral arguments:
\[ \tilde P_r(-n)=(-1)^{r-1}\frac{2}{n^2}\BC{2n}{n}^{-1}\sum_{k=1}^n\BC{2n}{n-k}k^{2r+1}.
\]
For example, one easily obtains $\tilde P_r(-1)=(-1)^{r-1}$ and $\tilde P_r(-2)=(-1)^{r-1}(2^{2r-1}+1)/3$.

Likewise, one can use the relation of the moments of the absolute distance to the origin in a symmetric Bernoulli random walk  and the sums $S_r(n)$ to express these moments in terms of the polynomials $P_r(n)$ and $Q_r(n)$:
\[ \Exp\abs{Y_{2n}}^{2r}=2^r\, Q_r(n)
   \ \ \ \mbox{\rm and} \ \ \ 
   \Exp\abs{Y_{2n}}^{2r+1}= \BC{2n}{n}2^{2(r-n)+1} 
     \, 
    n\,P_{r}(n). 
\]
This equivalence can be used to establish the rate of convergence to the moments of the half-normal distribution. 

Finally, it should be noted that one can also determine expressions for $S_{2r}(n)$ 
by means of the generating function
\[ f_n(\varphi) 
   =\sum_{k=0}^{2n}\BC{2n}{k}e^{(n-k)\varphi}
   =e^{n\varphi}\left[1+e^{-\varphi}\right]^{2n} 
   =2^n\left[1+\cosh \varphi \right]^n,
\]
so that $S_{2r}(n)=f_n^{(2r)}(0)$.
However, this approach covers only the even indexed case, 
and does not give the same insight into the problem as the one that we have followed here.

\section{Acknowledgment}
I would like to thank the anonymous referee for drawing attention to the occurrence of the Genocchi numbers
in the polynomials $P_r(n)$. This led to a further investigation and the characterization in terms of the Gandhi
and Dumont-Foata polynomials.

\par\noindent
AMS Classification Numbers: 
   11B65, 60G50, 44A60

\begin{thebibliography}{10}

\bibitem{Bruckman:1999}
Paul~S. Bruckman.
\newblock Problem {B-871}.
\newblock {\em The Fibonacci Quarterly}, 37(1):85, February 1999.

\bibitem{Carlitz:1972}
L.~Carlitz.
\newblock A conjecture concerning {G}enocchi numbers.
\newblock {\em Det Kongelige Norske Videnskabers Selskabs Skrifter}, 9:1--4,
  1972.
 
\bibitem{Carlitz:1980}
L.~Carlitz.
\newblock Explicit formulas for the {D}umont-{F}oata polynomial.
\newblock {\em Discrete Mathematics}, 30(3):211--225, 1980.
 
\bibitem{Comtet:1974}
Louis Comtet.
\newblock {\em Advanced combinatorics; the art of finite and infinite
  expansions}.
\newblock D. Reidel, Boston, 1974.

\bibitem{Dumont:1974}
Dominique Dumont.
\newblock Interpretations combinatoires des nombres de {G}enocchi.
\newblock {\em Duke Mathematical Journal}, 41:305--318, 1974.

\bibitem{DumontFoata:1976}
Dominique Dumont and Dominique Foata.
\newblock Une propri{\'e}t{\'e} de sym{\'e}trie des nombres de {G}enocchi.
\newblock {\em Bulletin de la Soci{\'e}t{\'e} Math{\'e}matique de France},
  104(4):433--451, 1976.

\bibitem{DynkinUspenskii:1963}
E.~B. Dynkin and V.~A. Uspenskii.
\newblock {\em Random Walks: Part Three of Mathematical Conversations}.
\newblock D. C. Heath and Company, Boston, 1963.

\bibitem{Feller:1950}
William Feller.
\newblock {\em An Introduction to Probability Theory and Its Applications}.
\newblock John Wiley \& Sons, Inc., New York, 3rd edition, 1950.

\bibitem{Gandhi:1970}
J.~M. Gandhi.
\newblock A conjectured representation of {G}enocchi numbers.
\newblock {\em The American Mathematical Monthly}, 77(5):505--506, May 1970.

\bibitem{GrahamKP:1989}
Ronald~L. Graham, Donald~E. Knuth, and Oren Patashnik.
\newblock {\em Concrete Mathematics}.
\newblock Addison-Wesley, Reading, Massachusetts, 1989.

\bibitem{Gut:1988}
Allan Gut.
\newblock {\em Stopped Random Walks}.
\newblock Springer-Verlag, New York, 1988.

\bibitem{JohnsonKB:1994}
Norman~L. Johnson, Samuel Kotz, and N.~Balakrishnan.
\newblock {\em Continuous Univariate Distributions}, volume~I.
\newblock John Wiley \& Sons, Inc., New York, 1994.

\bibitem{Revesz:1990}
P\'al R\'ev\'esz.
\newblock {\em Random Walk in Random and Non-Random Environments}.
\newblock World Scientific Publishing Co., Singapore, 1990.

\bibitem{RiordanStein:1973}
John Riordan and Paul~R. Stein.
\newblock Proof of a conjecture on {G}enocchi numbers.
\newblock {\em Discrete mathematics}, 5(4):381--388, 1973.

\bibitem{Sloane:2000}
N.~J.~A. Sloane.
\newblock {\em {The On-Line Encyclopedia of Integer Sequences}}, 2000.
\newblock Published electronically at
  http://www.research.att.com/$\sim$njas/sequences/.

\bibitem{Strazdins:2000}
Indulis Strazdins.
\newblock Solution to problem {B-871}.
\newblock {\em The Fibonacci Quarterly}, 38(1):86--87, February 2000.

\bibitem{Strehl:1979}
Volker Strehl.
\newblock Alternating permutations and modified {G}handi-polynomials.
\newblock {\em Discrete Mathematics}, 28:89--100, 1979.

\end{thebibliography}
\end{document}